\documentclass[12pt]{article}
\usepackage{epic,eepic,amsmath,amssymb,amsthm,fancybox,setspace}
\usepackage{latexcad}
       \newtheorem{theorem}{Theorem}
       
       \newtheorem{lemma}[theorem]{Lemma}

       {\theoremstyle{definition}
       \newtheorem{definition}{Definition}
       
       }
\def\N{\ensuremath{\mathbb N}} 
\def\T{\ensuremath{\mathbb T}} 
\def\R{\ensuremath{\mathbb R}} 
\def\Z{\ensuremath{\mathbb Z}} 

\newcommand{\ds}{\displaystyle}

\newcommand{\eqq}{\equiv}

\newcommand{\GGG}{\Gamma}
\newcommand{\hGGG}{{\widehat{\Gamma}}}

\newcommand{\LLL}{\Lambda}

\newcommand{\LOG}{\log}

\newcommand{\oo}{\infty}
\newcommand{\sse}{\subset}
\newcommand{\ess}{\emptyset}
\newcommand{\sm}{\setminus}

\newcommand{\I}{{\mathbf{1\!\!\!1}}}
\newcommand{\bT}{{\mathbf{T}}}
\newcommand{\bN}{{\mathbf{N}}}
\title{An $L^1$ Counting Problem
in Ergodic Theory}
\author{
Idris Assani,
 Department of Mathematics,\\
 University of North Carolina at Chapel Hill,\\
 Chapel Hill, North Carolina 27599, USA\\
email: assani@email.unc.edu\\
{\tt www.math.unc.edu/Faculty/assani}
\\
  \\
\\
Zolt\'an Buczolich\thanks{This work was completed while this
author visited the Departments of Mathematics of University of
North Texas and of University of North Carolina at Chapel
Hill.},
Department of Analysis, E\"otv\"os Lor\'and\\
University, P\'azm\'any P\'eter S\'et\'any 1/c,
1117 Budapest, Hungary\\
email: buczo@cs.elte.hu\\
{\tt www.cs.elte.hu/\hbox{$\sim$}buczo}\\
 \\
and\\
 \\
R. Daniel Mauldin\thanks{Supported in part by NSF grant DMS 0100078
},
 Department of Mathematics,\\
 University of North Texas,
 Denton, Texas 76203-1430, USA\\
email: mauldin@unt.edu\\
{\tt www.math.unt.edu/\hbox{$\sim$}mauldin}
}
\date{\today}
\begin{document}
\maketitle
\begin{abstract}
We solve the following counting problem for measure preserving
transformations. For $f\in L_+^1(\mu)$, is it true that $\ds
\sup_n\frac{\bN_n(f)(x)}{n} <\infty,$ where
$$\ds\bN_n(f)(x)=
\# \left \{k: \frac{f(T^k x)}{k}>\frac 1 n \right\}?$$ One of the
consequences is the nonvalidity of J. Bourgain's Return Time Theorem
for pairs of $(L^1, L^1)$ functions.
\end{abstract}

\doublespacing

\section{Introduction}
Let $(X, \mathcal{B}, \mu)$ be a probability measure space, $T$ an
 invertible measure preserving transformation on this space and
$f\in L_+^1(\mu)$. Since $\frac{f(T^nx)}{n}\rightarrow 0$ a.e., the
following function
$$\bN_n(f)(x)= \# \left \{k: \frac{f(T^k x)}{k}>\frac 1 n
\right\}$$ is finite
a.e. In this paper we consider the
following\\
 \vskip1ex
\noindent{\bf Counting Problem I.}
{\em Given $f\in L_+^1(\mu)$ do we have
$\sup_n\frac{\bN_n(f)(x)}{n}
 <\infty$,\\ $\mu$ a.e.?}\\
 \vskip1ex

 In \cite{[A1]} and \cite{[A2]} the operator $\sup_n\frac{\bN_n(f)(x)}{n}$
 was introduced and the pointwise convergence of
 $\frac{\bN_n(f)(x)}{n}$ was studied. It was shown there that if
 $f\in
L_+^p$ for $p>1$, or $f\in L\log L$ and the transformation $T$ is ergodic,
then $\frac{\bN_n(f)(x)}{n}$ converges a.e to $\int fd\mu$. If $T$
is not ergodic, then the limit is the conditional expectation of
the function $f$ with respect to the $\sigma$ field of the
invariant sets for $T$. Hence, the limit is the same as the limit
of the ergodic averages $\frac{1}{N}\sum_{n=1}^N f(T^nx)$. The
limit of the ergodic averages, by Birkhoff's pointwise ergodic
theorem, exists for any function $f\in L^1(\mu)$. It is natural
to ask whether $\frac{\bN_n(f)(x)}{n}$ also converges a.e.,
when $f\in L^1(\mu)$. Another motivation for this question is given by the fact that
for i.i.d. random variables $X_n\in L^1$ it was shown in
\cite{[A1]} that $$\frac{\#\{k:
\frac{X_k(\omega)}{k}>\frac{1}{n}\}}{n}$$ converges a.e. to
$E(X_1)$. The counting problem was afterwards discussed in
\cite{[JRW]}.
 \vskip1ex

One can see by using the methods of
\cite{[A1]}, for instance,
 that the convergence for all functions
$f\in L_+^1(\mu)$ will be guaranteed if one can answer the
following equivalent problem.\\
 \vskip1ex
 \noindent{\bf Counting Problem II.}
{\it Does there exist a finite positive
constant $C$ such that for all measure preserving systems and all
$\lambda >0$
$$\mu\left \{x: \sup_n\frac{\bN_n(f)(x)}{n}>\lambda
\right\}\leq
\frac{C}{\lambda}\|f\|_1?$$}\\
\vskip1ex
\noindent
Our main result will be to show that this equivalent
problem has a negative answer. More precisely we have
\begin{theorem}\label{th1}
$$\sup_{(X,\mathcal{B}, \mu,
T)}\sup_{\|f\|_1=1}\sup_{\lambda>0}\lambda\cdot \mu
\left \{x:
\sup_n\frac{\bN_n(f)(x)}{n}>\lambda\right\}= \infty.$$
\end{theorem}

This theorem answers then the question raised in \cite{[A1]}.
\vskip1ex We will also derive answers to some related problems. The
first consequence, linked to the study of
the maximal function
$\bN^{*}(f)(x) =
\sup_n\frac{\bN_n(f)(x)}{n},$ is what we call the return times for
the tail (of the Cesaro averages).
\begin{definition}
Let $(X, \mathcal{B},\mu, T)$ be a measure preserving system. The
Return Times for the Tail
Property
holds in $L^r(\mu)$, $1\leq r \leq
{\infty}$ if for each $f\in L^r(\mu)$  we can find a set $X_f$ of
full measure such that for all $x\in X_f$ for all measure
preserving systems $(Y, \mathcal{G}, \nu, S)$ and each $g\in
L^1(\nu)$ the sequence $\frac{f(T^nx)\cdot
g(S^ny)}{n} $ converges to
zero for a.e. $y$.
\end{definition}

A first consequence of Theorem \ref{th1} will be the following
\begin{theorem}\label{th2}
The Return Times for the Tail
Property
does not hold for $p=1$.
\end{theorem}

We observe that in \cite{[A1]} and \cite{[A2]} it was shown
that the Return Times for the Tail
Property
holds in $L^p$ for $1<p\leq
\infty$ and even in $L\log L$.
\vskip1ex A second consequence is a
solution to the $(L^1, L^1)$ problem mentioned in \cite{[A1]},
\cite{[A3]} and \cite{[Ru2]}. To explain this problem we need a
few definitions.
\begin{definition}
A sequence of scalars $a_n$ is said to be good universal for the
pointwise ergodic theorem (resp. norm convergence) in $L^r$,
$1\leq r \leq \infty$ if for all dynamical systems $(Y,
\mathcal{G}, \nu, S)$ the averages
$$\frac{1}{n}\sum_{k=1}^n a_k\cdot
g(S^ky)$$ converge  a.e. (resp. in
$L^r(\nu)$ norm).
\end{definition}
In \cite {[B1]}, 
\cite{[B15]}, and \cite{[B2]} J. Bourgain showed that given
$f\in L^{\infty}(\mu)$ the sequence $f(T^nx)$ is $\mu$ a.e. good
universal for the pointwise convergence in $L^1$.
Using H\"older's
inequality and the maximal inequality for the ergodic averages one
can extend his result to the pairs $(L^p, L^q)$ where $\frac{1}{p}
+ \frac{1}{q} =1$. This was mentioned in \cite{[Ru1]}. Bourgain's
Return Time Theorem strengthens Birkhoff's theorem on the product
space when the functions, $f$ and $g$, respect duality. That is,
if the function $f\in L^p(\mu)$ for some $1\leq p \leq \infty$,
then the set of convergence obtained from the Return Times
Theorem works
for all functions $g\in L^q(\nu)$, where $\frac{1}{p} +\frac{1}{q}
= 1$, hence
it is
a universal set. However, fixing $f$ and $g$, the
projection of the convergence set onto the first factor obtained
by Birkhoff's theorem depends on both functions.
A weakness of the
Return Time Theorem is that it does not address the case of $f\in
L^1$ and $g\in L^1$. Birkhoff's theorem, on the other hand,
guarantees convergence for $f\otimes g\in L^1\times L^1$,
$\mu\otimes \nu$-almost everywhere.

\vskip1ex In \cite{[A3]}
random stationary weights (i.i.d. random variables) were given for
which one could go ``beyond" the duality apparently imposed by the
use of H\"older's inequality.
It was also shown that given $f\in
L^1(\mu)$ the sequence $(f(T^nx))$ is $\mu$-a.e. good universal
for the $L^1$ norm. In \cite{[A1]} a Multiple Return Times
 Theorem
for $L^1$ i.i.d. random variables was obtained while in \cite{[Ru2]}
a Multiple Return Times theorem was proved for $L^{\infty}$
stationary processes.
The $(L^1,L^1)$ problem was the
following.\\
\vskip1ex
\noindent
{\bf $(L^1, L^1)$ Problem.} {\it Given $f \in
 L^1(\mu)$, is the sequence $(f(T^nx))$,  $\mu$-a.e. good universal
for the pointwise ergodic theorem in $L^1$?}\\
\vskip1ex
\noindent
A consequence of Theorem 2 is the following solution to the $(L^1,
L^1)$ problem

\begin{theorem}\label{th3}
Bourgain's Return Time Theorem does not hold for pairs of $(L^1,
 L^1)$ functions.
\end{theorem}

We also derive in Section \ref{sec4}
 some consequences in $L^1(\T)$
between the continuous analog of the maximal function
$\sup_n\frac{\bN_n(f)(x)}{n}$, namely

$$A(f)(x) = \sup_t t\cdot m\left \{ 0<y<x: \frac{f(x-y)}{y}>t
\right\},$$ or, analogously,

$$A(f)(x) = \sup_t t\cdot m
\left \{ 0<y<x: \frac{f(y)}{x-y}>t\right\},$$

\noindent and the one sided Hardy--Littlewood maximal function.

\section{Proof of Theorem \ref{th1}}\label{main}

In this section $\mu$ will denote Lebesgue measure
on $\R$ and $\log$ will denote logarithm in {\it base}
$2$.
An interval $I$ is a $2^{-R}$ grid interval if there is some
$j\in \Z$ such that $I=[j\cdot 2^{-R},(j+1)2^{-R}).$

\subsection{Basic systems}\label{bassec}

A ``life" function is a map
$\nu :\N\to \N$ such that for each $N\in \N$,
$\nu (N)> N.$ Given a life function $\nu$,
a gain constant $M > 3$, and a startup time $N_{1}$ we choose
 a
sequence $N_{2},...,N_{M}$ so that
\begin{equation}\label{*1}
{N_{l}}= {20} + \nu({N_{l-1}}),\  l=2,...,M.
\end{equation}.

Our aim in this section is to prove the following
\begin{lemma}\label{base} Suppose that
 a gain constant $M>3$, a life function $\nu,$
 a support constant
$S<2^M$,
and a startup time $N_{1}>\max\{10, M\}$
are given. Choose the sequence $N_2,...,N_{M}$ based
on $M$, $\nu$, and $N_{1}$ satisfying (1).
Given any $2^{-R}$ grid interval $I$,
there exists a positive integer
 $J_{0} > R$,
 disjoint subsets $\GGG_{1},...,\GGG_{M}$
of $I$,
and for each integer $J\geq J_{0}>R$
there is a  simple function $f:\R\to\R$,
 such that
$f(x)=0$ for $x\not\in I$ and
if
$T(x)=x+2^{-J}$ then  for all
$l=1,...,M$,
\begin{equation}\label{*base}
\frac{\bN_{n}(f)(x)}{n}>0.99\cdot 2^{-l+1}
\quad\text{ when }2^{N_{l}}\leq n\leq 2^{\nu(N_{l})}
\end{equation}
 for all $x\in \GGG_{l}$.
Moreover, each set
$\GGG_{l}$ consists of the union of intervals of the form
$[i\cdot 2^{-J_{0}},(i+1)2^{-J_{0}})$,
$\mu (\GGG_{l})>0.99\cdot 2^{-M+l-1}\mu (I)$, and
$\int_{I}f=2^{-M+1}\mu (I).$ We can also require that
$f(x)=0$ for any $x$ which is not in an interval
of the form $[(i\cdot 2^{M}+S)2^{-J}, (i\cdot 2^{M}+S+1)2^{-J})$
for some $i\in \Z.$
\end{lemma}

\begin{proof}

Set $h_{0}=2^{M+10}$ and choose $J_{0}$ such that
\begin{equation}\label{**1}
2^{10}2^{\nu({N_{M}})}h_{0}2^{-J_{0}}<2^{-R},
\end{equation}
or equivalently, $J_0 > \nu(N_M) + M + 20 + R.$
Now, let an integer $J\geq J_{0}$ be given.
Set $h=h_{0}\cdot 2^{J-J_{0}} = 2^{M+10+J-J_0}$. We shall first define
a sequence of sets $B_M, B_{M-1},\ldots,B_1$ each as the union of some
intervals in a corresponding sequence of finer dyadic grids. To begin
put
\begin{align}\label{*204}
B_{M}=&I\cap \bigcup_{j\in \Z}[2j\cdot 2^{-10-J}\cdot 2^{N_{M}}h
,(2j+1)\cdot 2^{-10-J}\cdot 2^{N_{M}}h )\\
\nonumber
=&I\cap \bigcup_{j\in \Z}[2j\cdot 2^{-10-J_{0}}\cdot 2^{N_{M}}h_{0}
,(2j+1)\cdot 2^{-10-J_{0}}\cdot 2^{N_{M}}h_{0} ).
\end{align}
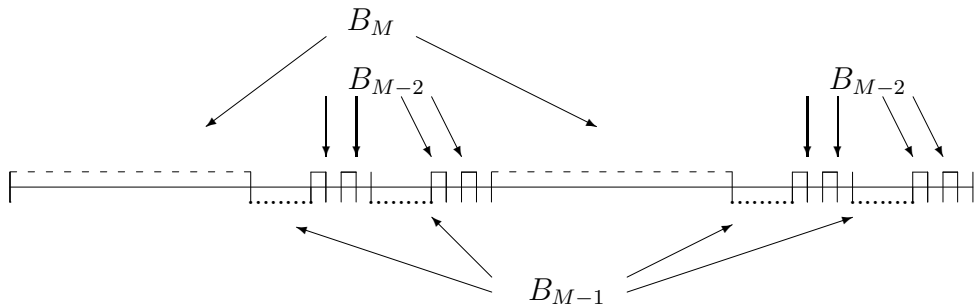
\placedrawing{ias1b.lp}{The sets $B_M$, $B_{M-1}$, and $B_{M-2}$}
{fig:BM}

Thus, $B_M$ consists of the intervals in the standard
$2^{N_M+M-J_0}$ grid with even index, $j$ and that are subsets of the
interval $I$. Clearly, $\mu (B_{M})=\mu (I)/2.$ In Figure
\ref{fig:BM} we illustrate the manner in which the sets $B_{l}$,
$l=M-2,M-1,M$, are located in $I$. Of course, in an illustration
we cannot divide an interval into several thousand pieces, so in
the figure the set $B_{M}$ consists of two intervals of length
$\mu(I)/4$, marked by dashed line, $B_{M-1}$ consists of four
intervals of length $\mu(I)/16$, marked by dotted line, $B_{M-2}$
consists of eight intervals marked by solid lines. The complement
of $B_M\cup B_{M-1}\cup B_{M-2}$ consists of eight ``unmarked"
intervals, each of the same length as the components of $B_{M-2}.$

In \eqref{*204}
the first expression for $B_M$ is given for some computational
purposes whereas the second expression shows that $B_M$ does not
depend on $J$ but rather on $J_0$. The same is true for all the
sets $B_i$ to be defined now.

Assume that $l\in  \{ 0,...,M-3\}$ and $B_{M-l'}$ is given for all
$l'\in  \{ 0,...,l\}.$
Set
\begin{align}\label{*206}
&B_{M-(l+1)}=\\ \nonumber &=(I\sm
\bigcup_{l'=0}^{l}
B_{M-l'})\cap \bigcup_{j\in \Z}[2j\cdot 2^{-10-J}\cdot 2^{N_{M-l-1}}h
,(2j+1)\cdot 2^{-10-J}\cdot 2^{N_{M-l-1}}h )=\\ \nonumber
&=(I\sm
\bigcup_{l'=0}^{l}
B_{M-l'})\cap \bigcup_{j\in \Z}[2j\cdot 2^{-10-J_0}\cdot 2^{N_{M-l-1}}h_0
,(2j+1)\cdot 2^{-10-J_0}\cdot 2^{N_{M-l-1}}h_0 )
.\end{align}

Thus, the set $B_{M-(l+1)}$ consists of the intervals with even
index in the standard $2^{N_{M-l-1}+M-J_0}$ grid that are
  subsets of $I$ and are not in $\cup_{i=M-l}^M B_i$.

Finally, if $B_{M-l}$ is given for $l\in  \{ 0,...,M-2  \}$,
we set $$B_{1}=B_{M-((M-2)+1)}=I\sm \bigcup_{l=0}^{M-2}B_{M-l}.$$
Returning to the illustration on Figure \ref{fig:BM}, if
$M=4$ then $B_M=B_4$ is marked by the dashed line,
$B_3$ is by the dotted line, $B_2$ by the solid line,
and $B_1$, the complement of the other three is the
``unmarked" part of $I$.

Observe that $\mu (B_{M-l})=\mu (I)/2^{l+1}$
holds for $l=0,...,M-2$ and $\mu (B_{1})=\mu (I)/2^{M-1}>
\mu (I)/2^{M}=\mu(I)/2^{(M-1)+1}.$
The set $B_{1}$ is the union of some disjoint intervals of the form
\begin{align}
[(2j-1)\cdot 2^{-10-J}\cdot 2^{N_{2}}h
&,2j\cdot 2^{-10-J}\cdot 2^{N_{2}}h )= \label{red6}\\
\nonumber =[(2j-1)\cdot 2^{-10-J_0}\cdot 2^{N_{2}}h_0
&,2j\cdot 2^{-10-J_0}\cdot 2^{N_{2}}h_0 ),
\end{align}
while for any $l=0,...,M-2$ the set $B_{M-l}$ is the union of some intervals
of the form
\begin{align}
[2j\cdot 2^{-10-J}\cdot 2^{N_{M-l}}h
&,(2j+1)\cdot 2^{-10-J}\cdot 2^{N_{M-l}}h )
\label{red7}\\ \nonumber
[2j\cdot 2^{-10-J_0}\cdot 2^{N_{M-l}}h_0
&,(2j+1)\cdot 2^{-10-J_0}\cdot 2^{N_{M-l}}h_0 ).
\end{align}
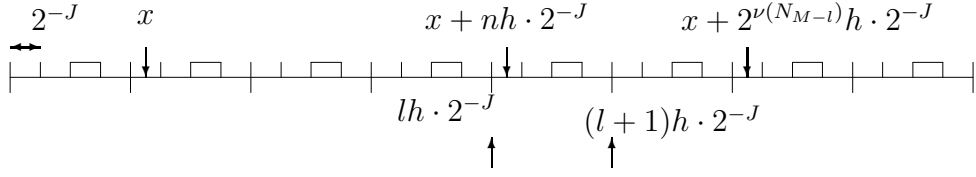
\placedrawing{ias2.lp}{The definition of $f$ in an interval $I'$}
{fig:ias2}

Our function $f$ which depends on $J$ will have value $0$ on $\cup_{l=2}^{M}B_{l}.$ To
determine its values on $B_1$, consider one of the intervals making up $B_1$:
$$I'=[(2j-1)\cdot 2^{-10-J}\cdot 2^{N_{2}}h
,2j\cdot 2^{-10-J}\cdot 2^{N_{2}}h )\sse B_{1}.$$
For each $l\in \Z$ such that the interval $[lh\cdot 2^{-J},
(l+1)h\cdot 2^{-J})\sse I'$ (and there are $\frac{2^{N_2}}{2^{10}}$ such $l$), choose exactly one $l'$
such that $lh\leq l'<
(l+1)h$, $l'\eqq S$ modulo $2^M$ and set $f(x)=h$
for $x\in [l'\cdot 2^{-J},(l'+1)\cdot 2^{-J})$,
otherwise we set $f(x)=0$.

In Figure \ref{fig:ias2} one can see one interval $I'$ being
enlarged. Again we could not divide this interval in a drawing
into several thousand subintervals, so in this illustration $h=4$,
and $S=2$. One tiny interval is of length $2^{-J}$, the tiny
intervals marked by an extra solid line are the ones where $f=h.$

From the definition of $f$,
we have
$\int_{I'}f= \frac{h}{2^J}\cdot\frac{2^{N_2}}{2^{10}}
 = \mu (I')$. By summing this
over all subintervals of $B_{1}$ of type $I'$,
we obtain $\int_I f=\int_{B_{1}}f=\mu (B_{1})=\mu (I)/2^{M-1}.$

Suppose $2^{N_{1}}\leq n\leq 2^{\nu({N_{1}})}$, and
\begin{equation}\label{*3}
[x,x+h\cdot 2^{\nu({N_{1}})-J})\sse I'.
\end{equation}
Then
$N_{1}>10$ implies $1000\leq n$ and hence
\begin{align*}
\bN_{n}(f)(x)
=& \#\{ k:\frac{f(T^{k}x)}{k}>\frac{1}{n}  \}=\\
=&
 \#\{ k:hn>k \text{ and }f(T^{k}x)=h  \}> 0.99\cdot \frac{nh}h=0.99n,
\end{align*}
Of course, instead of $0.99$ we could have used $0.999$, but this
is not of any consequence for our purposes.

Now, we define the
sets $\Gamma_i$ which do not depend on $J$ from the sets $B_i$. To begin set
\begin{align}\label{*203}
\GGG_{1}=& \{ x\in B_{1}: [x,x+h\cdot 2^{\nu({N_{1}})-J})
\sse B_{1}  \}=\\ \nonumber
=&\{ x\in B_{1}: [x,x+h_0\cdot 2^{\nu({N_{1}})-J_0})
\sse B_{1}  \}
.
\end{align}
Again, the second expression here shows that $\Gamma_1$ does not
depend on $J$ since $B_1$ does not depend on $J$.
For each interval $I'$ making up $B_1$, by using \eqref{*1}, we have
$|\Gamma_1 \cap I'| \geq |I'| - h\cdot 2^{\nu(N_1)-J} \geq
|I'|\cdot(1-2^{-(N_2-\nu(N_1)-10)}) > 0.99|I'|.$ So,
$\mu (\GGG_{1})>0.99 \cdot \mu (B_{1}).$

Observe that
for each $l=1,...,M-2$, the set $I\sm \cup_{i=0}^{l-1}
B_{M-i}=\cup_{i=1}^{M-l}B_{i}$ is the union of some intervals
of the form
$$I'_{M-l}=[(2j-1)\cdot 2^{-10-J}\cdot 2^{N_{M-l+1}}h
,2j\cdot 2^{-10-J}\cdot 2^{N_{M-l+1}}h ).$$

Also, the two sets $B_{M-l}$ and $B_1\cup\ldots\cup B_{M-l-1}$ are equally
distributed in $I_{M-l}'$ in the sense that if one takes the
$2^{N_{M-l}}h/2^{10+J}$ grid of the
interval $I_{M-l}'$, then every evenly indexed  interval is in $B_{M-l}$ and
the others are in $B_1\cup\ldots\cup B_{M-l-1}$. In particular,
$\mu (B_{M-l}\cap I_{M-l}')=\mu (I_{M-l}')/2
 =\mu (\cup_{i<M-l}B_{i}\cap I_{M-l}').$

Finally, by induction one can also see that
\begin{align}
\mu (B_{1}\cap I'_{M-l})=\mu (I'_{M-l})/2^{M-l-1},
\end{align}
and, more generally, if $n\in [2^{N_{M-l}},2^{\nu(N_{M-l})}]$
and $[x,x+nh\cdot 2^{-J})\subset I'_{M-l}$, then
\begin{align}\label{*201}
\mu (B_{1}\cap [x,x+nh\cdot 2^{-J}))
> 0.995nh\cdot 2^{-J}/2^{M-l-1}.
\end{align}
Set
\begin{align}\label{*2031}\GGG_{M-l}&= \{ x\in B_{M-l}:
[x,x+h\cdot 2^{\nu({N_{M-l}})-J})
\sse\bigcup_{l'\leq M-l} B_{l'}  \}
=\\
\nonumber
&=\{ x\in B_{M-l}: [x,x+h_{0}\cdot 2^{\nu({N_{M-l}})-J_{0}})
\sse\bigcup_{l'\leq M-l} B_{l'}  \}.
\end{align}
Using \eqref{*1}
 one can see that $\mu (\GGG_{M-l})>0.99 \mu (B_{M-l})\geq
0.99\mu(I)/2^{l+1}$.
If $x\in \GGG_{M-l}$ and $I_{M-l}'$ is the subinterval of $\cup_{l'=1}
^{M-l}B_{l'}$ containing $x$, then $x+jh\cdot 2^{-J}\in I'_{M-l}$
for all $0\leq j \leq 2^{\nu(N_{l})}.$
By using \eqref{*201} and the definition of $f(x)$ we have
 $$\bN_{n}(f)(x)= \#\{ k: hn>k \text{ and } f(T^{k}x)= h \}
\geq 0.99 \frac{nh}{h\cdot 2^{M-l-1}}=0.99 \frac{n}{2^{M-l-1}}.$$
From
$N_{2}>N_{1}>10$,
\eqref{*204}, \eqref{*206},
\eqref{red6}, \eqref{red7},
\eqref{*203}, and \eqref{*2031}
it follows that each $\GGG_{l}$ is the union of intervals
of the form $[i\cdot 2^{-J_{0}}, (i+1)\cdot 2^{-J_{0}}).$

\end{proof}

\subsection{Level $k$ systems}
\label{levelksec}

In this section the gain constant $M\in \N$ is fixed.

Next we define the life functions for all $k\in \N$.
We will use these functions in the proof of Lemma \ref{levelk}.
Set $\nu _{1}(N)=N+1$ for any $N\in \N.$
We proceed by induction, so assume that for $k\in \N$
we have already defined $\nu _{k}.$
If some $N\in \N$ is given use $\nu =\nu _{k}$
and $N_{1}=N_{1}^{(k)}(N)=N$ in \eqref{*1} to determine
the sequence $N_{2}^{(k)}(N),...,N_{M}^{(k)}(N).$
Put $\nu _{k+1}(N)=\nu _{k}(N_{M}^{(k)}(N))>N.$

We say that a random variable $X:I\to\R$ is
$(M-0.99)$-distributed
on $I$
 if $X(x)\in  \{0, 0.99,0.99\cdot 2^{-1},...,
0.99\cdot 2^{-M+1}  \}$
and $\mu ( \{ x:X(x)=0.99\cdot 2^{-l+1}  \})=0.99\cdot 2^{-M+l-1}
\mu (I),$
for $l=1,...,M.$

This section is about the existence of level $k$ systems, by which
we mean any system $(T,f)$ satisfying the conditions described in
the next lemma.

\begin{lemma}\label{levelk}
For any
$2^{-R}$ grid interval $I_{0}$,
positive integer $k\leq 2^{M}$, and any startup time
$K_{S}^{(k)}>\max \{ 10,M  \}$ there exists $J_{0}>0$ such that
for all $J\geq J_{0}$ we can find a system $(T,f)$ with the
following properties. The transformation $T$ is given by
$T(x)=x+2^{-J}.$ We have independent
$(M-0.99)$-distributed random variables $X_{h}$, $h=1,...,k$,
on $I_0$
 and an exit
time $K_{e}^{(k)}$
such that  for any $x\in I_{0}$
there exists an $n\in [2^{K_{S}^{(k)}},
2^{K_{e}^{(k)}}]$,
for which
\begin{equation}\label{ergosumest}
\frac{\bN_{n}(f)(x)}{n}\geq \sum_{h=1}^{k}X_{h}(x).
\end{equation}
 Moreover, $f$ is constant on the intervals of the form $[i\cdot 2^{-J},
(i+1)2^{-J})$, $\int_{I_{0}}f=k \cdot 2^{-M+1}\mu(I_0)$,
$f(x)=0=X_{h}(x)$
for $x\not\in I_{0}$, $h=1,...,k.$ We also may require that if
$$x \not \in  \bigcup_{l=0}^{k-1}\bigcup _{i\in\Z} [(i\cdot
2^{M}+l)2^{-J},(i\cdot 2^{M}+l+1)\cdot 2^{-J}),$$ then $f(x)=0.$
\end{lemma}

\begin{proof}
To define our level $1$ systems we use Lemma \ref{base} on
$I_{0}$.
We
apply Lemma \ref{base} with $\nu =\nu _{1}$,
and
$N_{1}=K_{S}^{(1)}$. So,
$K_{e}^{(1)}=\nu _{1}({N_{M}})$ will be the exit time. We choose
our $(M-0.99)$-distributed random variable the following way. For
$l=1,...,M$ we select a measurable set $\hGGG_{l}\sse \GGG_{l}$
such that $\mu (\hGGG_{l})=0.99 \cdot 2^{-M+l-1}\cdot 2^{-R}$. If $x\in
\hGGG_{l}$ for some $l$ then we set $X_{1}(x)=0.99\cdot 2^{-l+1}$
and $X_{1}(x)=0$ otherwise. Viewed in this way Lemma \ref{base}
guarantees that level $1$ systems exist.

We proceed by induction on $k$. Assume that level $k$ systems exist and we need to
verify the existence of level $k+1$ systems.

First, calling upon Lemma \ref{base}, we define a ``mother" base
system. The ``subsystems" of this ``mother" system will be level
$k$ systems with different life intervals. Here is a heuristic
argument behind our construction. Due to the $L^{1}$ restrictions,
the mother system is unable to deal with all the subsystems
simultaneously at the same time. So some subsystems have longer
and longer waiting times, but the longer the waiting time, the
longer lifetime they need. Since we already know how the
subsystems will look, this information is encoded by the life
function $\nu_{k+1}$. Now, using this function, we can ``design" a
mother system which can accomodate all the subsystems. Let us
proceed.

Given the startup constant $N_{1}=N_{1,0}=K_{S}^{(k+1)}>\max  \{ 10,M  \}$
putting the life function $\nu_{k+1}$
defined at the beginning of Subsection \ref{levelksec}
into \eqref{*1}, determine the sequence
$N_{2,0},...,N_{M,0}$, (the extra $0$ in subscripts will refer to
the ``mother system"). We
also put $N_{0,0}=N_1,$ and
 set the support constant $S_{0}=k$
for the mother system.

Next we apply Lemma \ref{base}
with $\nu =\nu _{k+1}$
to the $2^{-R}$ grid interval
$I_{0}=[j_{0}\cdot 2^{-R},(j_{0}+1)\cdot 2^{-R})$ we choose $J_{0,0}$ and
 disjoint subsets $\GGG_{1,0},...,\GGG_{M,0}$
of $I_0$ such that for each $l=1,...,M$, $\GGG_{l,0}$ consists of
the union of some intervals of the form $[i\cdot
2^{-J_{0,0}},(i+1)2^{-J_{0,0}})$, and $\mu (\GGG_{l,0})>0.99\cdot
2^{-M+l-1}\cdot 2^{-R}.$ For any $J\geq J_{0,0}$ we
can choose a function $\phi
_{0}=f:I_{0}\to\R$,
 such that
 if $T(x)=x+2^{-J}$ then  for all $l=1,...,M$,
\begin{equation}\label{*base0}
\frac{\bN_{n}(\phi _{0})(x)}{n}>0.99\cdot 2^{-l+1},
\quad\text{ when }2^{N_{l,0}}\leq n\leq 2^{\nu_{k+1}(N_{l,0})},
\end{equation}
 for all $x\in \GGG_{l,0}$.
Moreover, $\int_{I_0}\phi_{0}=2^{-M+1}\mu(I_0).$ Since $S_{0}=k$, we also
have $\phi _{0}(x)=0$ for any $x$ which is not in an interval of
the form $[(i\cdot 2^{M}+k)2^{-J}, ((i\cdot 2^{M}+k+1)2^{-J})$ for
some $i\in \Z.$

Next, consider the intervals $I_{j}=[j_{0}\cdot 2^{-R}+ (j-1)\cdot
2^{-J_{0,0}},j_{0}\cdot 2^{-R}
+j\cdot 2^{-J_{0,0}} )$ for $j=1,...,2^{J_{0,0}-R}$. Our
``subsystems" will live on these intervals.

If $I_{j}\sse \cup_{l=1}^{M}\GGG_{l,0}$
then
there is a unique $l(j)$
such that $I_{j}\sse \GGG_{l(j),0}.$
If $I_{j}\not\sse \cup_{l=1}^{M}\GGG_{l,0}$
then
$I_{j}\cap \cup_{l=1}^{M}\GGG_{l,0}=\ess$, and
in this case we set $l(j)=0.$
By our assumption
on any $I_{j}$
we can find
level $k$ systems. So, for each $j\in  \{ 1,...,2^{J_{0,0}-R}  \}$
choose a level $k$ system on $I_{j}$ with startup time
$K_{S,j}^{(k)}=N_{l(j),0}.$ Choose $J_{0,j}$ for each $j=1,...,
2^{J_{0,0}-R}$ according to our induction hypothesis. Set
$J_{0}=\max{ \{ J_{0,j}:j=0,...,2^{J_{0,0}-R}  \}}$ and choose a
$J\geq J_{0}$. The transformation $T$ will be given by
$T(x)=x+2^{-J}.$ For this $J$ choose $\phi _{0}$ as was explained
above, and by the induction hypothesis for any
$j=1,...,2^{J_{0,0}-R}$ choose $\phi _{j}=f$ and independent
$(M-0.99)$-distributed random variables $X_{h,j}$, $h=1,...,k$,
on $I_j$,
 and an exit
time $K_{e,j}^{(k)}=\nu_{k+1}(N_{l(j),0})$
such that  for any $x\in I_{j}$
there exists an $n\in [2^{N_{l(j),0}},
2^{\nu_{k+1}(N_{l(j),0})}]$,
for which
\begin{equation}
\frac{\bN_{n}(\phi _{j})(x)}{n}\geq \sum_{h=1}^{k}X_{h,j}(x).
\end{equation}
 Moreover, $\phi _{j}$ is constant on the
 intervals of the form $[i\cdot 2^{-J},
(i+1)2^{-J})$, $\int_{I_{j}}\phi _{j}=k \cdot 2^{-M+1}\mu
(I_{j})$, $\phi _{j}(x)=0=X_{h,j}(x),$ for $x\not\in I_{j}$,
$h=1,...,k.$ We may also require that if $$x \not \in
\bigcup_{l=0}^{k-1}\bigcup _{i\in\Z} [(i\cdot
2^{M}+l)2^{-J},(i\cdot 2^{M}+l+1)\cdot 2^{-J})$$ then $\phi
_{j}(x)=0.$ This last property implies that the support of
$\phi_{0}$ is disjoint from the support of any $\phi _{j},$
$j=1,...,2^{J_{0,0}-R}$. Since $\phi _{j}$ is supported on $I_{j}$,
we see that the supports of the functions $\phi _{j}$ are also
disjoint.

Set $f=\sum_{j=0}^{2^{J_{0,0}-R}}\phi _{j}.$ Then, using the fact
that the supports are disjoint, we have
$\bN_{n}(f)(x)=\sum_{j=0}^{2^{J_{0,0}-R}}\bN_{n}(\phi _{j})(x).$ We also
calculate
\begin{align*}
&\int_{I_{0}}f=\int_{I_{0}}\phi _{0}+ \sum_{j=1}^{2^{J_{0,0}-R}}
\int_{I_{j}}\phi _{j}=\\ &=2^{-M+1}\mu(I_{0})+ k\cdot 2^{-M+1}
\sum_{j=1}^{2^{J_{0,0}-R}} \mu (I_{j})= (k+1) 2^{-M+1}\mu(I_{0}).
\end{align*}

For $h=1,...,k,$ set $X_h=X_{h}'=\sum_{j=1}^{2^{J_{0,0}-R}}X_{h,j}$.
Let $X_{k+1}'(x)=0.99\cdot 2^{-l+1}$ if $x\in \GGG_{l,0}$,
otherwise set $X_{k+1}'=0.$ Since $X_{k+1}'$ is constant on the
intervals $I_{j},$ one can also see that the functions
$X_{h}'(x)$, $h=1,...,k+1$ are independent. The functions
$X_{h}'(x)$ are $(M-0.99)$-distributed on $I_{0}$ for $h=1,...,k$.
The function $X_{k+1}'(x)$ is not $(M-0.99)$-distributed, but is
$(M-0.99)$-superdistributed. By this we mean that $\mu ( \{
x:X_{k+1}'(x)=0.99\cdot  2^{-l+1}  \})\geq 0.99 \cdot 2^{-M+l-1}
\mu(I_{0}),$
for any $l=1,...,M.$ But we can and do choose $X_{k+1}\leq
X_{k+1}'$ such that $X_{k+1}$ is $(M-0.99)$-distributed on $I_{0}$
and the system $X_{h}(x)$, $h=1,...,k+1$ is independent.

If $x\in I_{j}\sse \GGG_{l(j),0},$ then
$$\frac{\bN_{n}(\phi _{0})(x)}{n}>0.99\cdot 2^{-l(j)+1}=X_{k+1}'(x)
\geq X_{k+1}(x),$$
when $2^{N_{l{(j),0}}}\leq n\leq 2^{\nu_{k+1}(N_{l(j),0})}.$
For these same $x$, by our induction hypothesis, there exists
$n\in [2^{N_{l(j),0}},2^{\nu_{k+1}(N_{l(j),0})}]$ for which
$$\frac{\bN_{n}(\phi _{j})(x)}{n}\geq \sum _{h=1}^{k}X_{h,j}(x)=
\sum_{h=1}^{k}X_{h}(x).$$
Therefore, there exists
$n\in [2^{N_{l(j),0}},2^{\nu_{k+1}(N_{l(j),0})}]
\sse [2^{K_{S}^{(k+1)}},2^{\nu _{k+1}(N_{M,0})}]$ for which
$$\frac{\bN_{n}(f)(x)}{n}=
\sum_{j=0}^{2^{J_{0,0}-R}}\frac{\bN_{n}(\phi _{j})(x)}{n} \geq
\sum_{h=1}^{k+1}X_{h}(x).$$
This also shows that the exit time $K_{e}^{(k+1)}$ can be
chosen to be $\nu _{k+1}(N_{M,0}).$

\end{proof}

\subsection{$p$-blocks}\label{pblocksec}

We restate in our measure theoretical language formula (9) on
p. 21 of \cite{[L]} in the form of a lemma.

\begin{lemma}\label{lamperti}
Assume that for a given $q\in \N$ we have independent identically
distributed random variables $X_{1},...,X_{q}$
on a probability space $(\Omega ,\Sigma,\mu )$, each
with finite mean $u$ and variance $v$.
Then for each $\epsilon >0$ we have
\begin{equation}\label{lampertieq}
\mu \left ( \left \{ x: \left |\left (\sum_{h=1}^{q}X_{h}(x)\right)-qu
\right|\geq q\epsilon  \right \}\right)\leq
\frac{q v}{(q\epsilon )^{2}}.
\end{equation}

\end{lemma}

This section concerns the existence of $p$-blocks as described in the
next lemma. We assume $I=[0,1).$

\begin{lemma}\label{pblock}
There exists $p_{0}>2$ such that for every
$p>p_{0}$ we can choose a $p$-block.
By this we mean, that we can find
a system $(T_{p},f_{p})$,
such that $
\frac1{p\LOG^{2}(p)}\leq
\int_{I}f_{p}\leq
\frac{4}{p\cdot \LOG^{2}(p)},$
 $T_{p}(x)=x+2^{-J_p}$, $\text{mod } 1$ for a large integer $J_p$.
There is a set $\LLL_{p}$ with $\mu (\LLL_{p})> 0.99$
and
there exists an exit time $E_{p}>2^{p}$ such that
for each $x\in \LLL_{p}$ there is some
 $n\in[2^{2^{p}},2^{E_{p}}]$ for which
$$
\frac{\bN_{n}(f_{p})(x)}{n}
>\frac{1}{2^2\LOG^{2}(p)}.$$
\end{lemma}

\begin{proof}
Using Lemma $ \ref{levelk}$
on $I_0=I=[0,1)$
 with $k=2^{p}$,
$M=M_{p}=[p+\LOG(p)+\LOG(\LOG^{2}(p))]$, and startup time
$K_{S}^{(2^{p})}=2^{p}$ we choose and fix $J_p\geq J_{0}$ and
 a level $2^{p}$
system $(T_{p},f _{p})$ such that $T_{p}(x)=
x+2^{-J_p}$, $\text{mod } 1$.
Here we remark that Lemma \ref{levelk} uses
$T_{p}(x)=x+2^{-J_{p}}$, but $f_{p}$ is supported
on $I$ and hence
by using
$T_{p}(x)=x+2^{-J_{p}}$,
$\text{mod } 1$
we cannot decrease $\bN_{n}(f_{p}).$
We have independent
$(M_{p}-0.99)$-distributed
random variables $X_{h,p},$ $h=1,...,2^{p}$,
 such that
for any $x\in I_{0}$ there exists $n\in [2^{2^{p}},2^{K_{e}^{(2^{p})}}]$
for which
$$\frac{\bN_{n}(f_{p})(x)}
{n}
\geq \sum_{h=1}^{2^{p}}X_{h,p}(x),$$
and $\int_{I} f_{p}=2^{p}\cdot 2^{-M_{p}+1}.$

Then, for any $h$,
\begin{align}\label{*CD3}
u=\int_{I}X_{h,p}= & \sum_{l=1}^{M_{p}}
0.99\cdot 2^{-l+1}\cdot 0.99\cdot  2^{-M_{p}+l-1}\\
\nonumber
\geq &\frac{0.99^{2}}{2^{p}p\cdot \LOG^{2}(p)}\cdot p >\frac1{2^{p+1}
\LOG^{2}(p)}.
\end{align}
and
\begin{align*}
v_{0}=&\int_{I}X_{h,p}^{2}= \sum_{l=1}^{M_{p}}
0.99^{2}\cdot 2^{-2l+2}
\cdot 0.99\cdot  2^{-M_{p}+l-1}=\\
&\sum_{l=1}^{M_{p}} 0.99^{3}\cdot 2^{-M_{p}-l+1}<
0.99^{3}\cdot  2^{-M_{p}+1}\sum_{l=1}^{\oo}2^{-l}\leq
\frac{4}{2^{p}\cdot p\cdot \LOG
^{2}(p)}.
\end{align*}
We have
$$
0<v=\int_{I}(X_{h,p}(x)-u)^{2}dx=v_{0}-u^{2}<v_{0}\leq
 \frac{4}{2^{p}\cdot p
\cdot \LOG^{2}(p)}.$$

Next, by Lemma \ref{lamperti} used with $\epsilon =1/2^{p+2}\LOG^{2}(p),$
$q=2^{p}$ we obtain that
\begin{align*}
&\mu \left ( \left \{ x: \left |\left (
\sum_{h=1}^{2^{p}}X_{h,p}(x)
\right)-
2^{p}u
\right|\geq 2^{p}\cdot
\frac{1}
{2^{p+2}\LOG^{2}(p)}
  \right \}\right)\leq
\\&
\frac{2^{p}\cdot \frac{4}{2^{p}\cdot p\cdot \LOG^{2}(p)}}
{(2^{p}\frac{1}
{2^{p+2}\cdot \LOG^{2}(p)})^{2}}=
\frac{64\cdot \LOG^{2}(p)}{p}.
\end{align*}
By using \eqref{*CD3} this implies
\begin{align*}
\mu \left ( \left \{ x:
\sum_{h=1}^{2^{p}}X_{h,p}(x)\leq
\frac{1}{2^{2}\LOG^{2}(p)}  \right \}\right)\leq
\frac{64\cdot \LOG^{2}(p)}{p}.
\end{align*}
Assume that $p$ is chosen so large that $64\cdot \LOG^{2}(p)/p<0.01.$
Then letting $$\LLL_{p}= \left \{ x\in I_{0}:\sum_{h=1}^{2^{p}}
X_{h,p}(x)>\frac{1}{2^{2}\LOG^{2}(p)}  \right\},$$
we have $\mu (\LLL_{p})>0.99.$
We set $E_{p}=K_{e}^{(2^{p})}$.
For any $x\in \LLL_{p}$ we have an $n\in [2^{2^{p}}, 2^{E_{p}}]$
such that
$$\frac{\bN_{n}(f_{p})(x)}{n}
\geq \sum_{h=1}^{2^{p}}X_{h,p}(x)>
\frac{1}{2^{2}\LOG^{2}(p)},$$
and
$$\frac{1}
{p\LOG^{2}(p)}
\leq\int_{I}f_{p}=
2^{-[p+\LOG(p)+\LOG(\LOG^{2}(p))]+1}\cdot 2^{p}\leq
\frac{4}{p\LOG^{2}(p)}.$$

\end{proof}

Next we turn to the proof of Theorem \ref{th1}.

\begin{proof}
By using Lemma \ref{pblock}
with $I=X=[0,1)$
choose $p_{0}$ and for each
$p>p_{0}$ a $p$-block.
Set $\phi _{p}=f_{p}/\int _{I}f_{p},$
and $\lambda _{p}=1/(2^{3}\LOG^{2}(p)\cdot \int_{I}f_{p})\geq
p/32$. By Lemma \ref{pblock}
for any $x\in\LLL_{p}$ there is $n'\in [2^{2^{p}},2^{E_{p}}]$
such that
$$\frac{\bN_{n'}(f_{p})(x)}{n'}>\frac{1}{2^{2}\LOG^{2}(p)}.$$
Now using the definition of $\bN_{n'}$
we obtain
\begin{align*}
&\frac{\bN_{n'}(f_{p})(x)}{n'}=
\frac{\# \{ k:f_{p}(T^{k}x)/k>1/n'  \}}{n'}=
\frac{\# \{ k:\phi_{p}(T^{k}x)/k>1/(n' \int_{I}f_{p}) \}}{n'}<\\
&
\frac{\# \{ k:\phi_{p}(T^{k}x)/k>1/([n' \int_{I}f_{p}]+1) \}}{[n'
\int_{I}f_{p}]+1}
\cdot \frac{n'\int_{I}f_{p}+1}{n'}=\\
\intertext{using $n=[n' \int_{I}f_{p}]+1$}
\\
&\frac{\bN_{n}(\phi_{p})(x)}{n}\left (\int_{I}f_{p}+\frac{1}{n'}\right)<
\frac{\bN_{n}(\phi_{p})(x)}{n}\cdot 2\int_{I}f_{p}.
\end{align*}
Hence for all $x$ from a set of measure at least $0.99$
there exist $n$ such that
$\bN_{n}(\phi _{p})(x)/n>\frac{1}{2^{3}\LOG^{2}(p)\int_{I}f_{p}}=\lambda _{p}.$

Since $\lambda _{p}\cdot 0.99\to\oo$ we have established
Theorem \ref{th1}.
\end{proof}

\section{Proofs of Theorems \ref{th2} and \ref{th3}}\label{sec3}

\subsection{Proof of Theorem \ref{th2}}\label{sec31}
\begin{proof}
 Theorem \ref{th2} follows from Theorem 8 in \cite{[A1]}. It was shown
 there that for a sequence of nonnegative numbers $c_n$ such that
 $\displaystyle \lim_n {c_n}/{n}=0$
the following two statements are equivalent
\begin{enumerate}
\item $$\sup_n\frac{\#\{k: \frac{c_k}{k}> \frac{1}{n}\}}{n}<\infty;
$$ and
\item for all measure preserving systems $(Y, \mathcal{G}, \nu,
S)$ and all $g\in L^1(\nu)$, the sequence ${c_n\cdot g(S^ny)}/{n}$
converges to zero $\nu$ a.e.
\end{enumerate}
Taking the sequence $c_n =f(T^nx)$ for an ergodic transformation
$T$  shows that if the validity of the Return Time for the
Tail Property in $L^1$ were to hold,
then we should have for all $f\in L_+^1(\mu)$ for a.e. $x$,
\begin{equation}
\label{idr}
\sup_n\frac{\#\{k: \frac{f(T^kx)}{k}> \frac{1}{n}\}}{n}<\infty.
\end{equation}
Condition \eqref{idr} for all $f\in L_+^1(\mu)$  is
equivalent to saying that
      $$\sup_{\alpha>0} \frac{\#\{k: \frac{f(T^kx)}{k}> \frac{1}
{\alpha}\}}{\alpha}<\infty$$ for all
 $f\in L_+^1(\mu)$ for a.e. $x$.
Consider an enumeration of the positive rational numbers $r_k$ and
define for each $k$ the function $\bT_k(f)(x) =
\frac{\bN_{r_k}(f)(x)}{r_k}$. We have
$$\sup_{\alpha>0} \frac{\#\{k: \frac{f(T^kx)}{k}> \frac{1}{\alpha}\}}
{\alpha} = \sup_k\bT_k(f)(x)$$
 When $T$ is ergodic it commutes with the family of powers of $T$. By
the ergodic theorem this family is mixing. Indeed,
 we have $$\lim_N \frac{1}{N}\sum_{n=1}^N \mu
(A\cap T^{_n}(B)) = \mu(A)\mu(B)$$ so for each
 $\rho\geq 1$ there exists a $n$ such that $\mu(A\cap T^{_n}(B))\leq
\rho \mu(A)\mu(B)$.
For each $\gamma \geq 1$ we have $$\sup_k\bT_k(\gamma f)(x) = \gamma
\sup_k\bT_k(f)(x).$$ Thus we can apply Theorem 4 of \cite{[A1]}
 to
conclude that there exists a finite positive constant $C$ such
that for all $f\in L^1_+$,
$$\mu\{x: \sup_k \bT_k(f)(x)>1\}\leq C\int f d\mu.$$  This means that
$$\mu\left \{x:\sup_{\alpha>0} \frac{\#\{k: \frac{f(T^kx)}{k}> \frac{1}
{\alpha}\}}{\alpha}> 1\right\}\leq C\int fd\mu.$$
Replacing the function $f$ by $f/{\lambda}$ provides a maximal
inequality for the maximal function
$$\sup_{\alpha>0} \frac{\#\{k: \frac{f(T^kx)}{k}> \frac{1}{\alpha}\}}
{\alpha}.$$
From this we obtain easily a maximal inequality with the same constant
$C$
for
$$\sup_n\frac{\#\{k: \frac{f(T^kx)}{k}> \frac{1}{n}\}}{n}.$$
Having this constant for one ergodic transformation provides the same
constant for all ergodic transformations.
 The ergodic decomposition would then show that
 $$\sup_{(X,\mathcal{B}, \mu,
T)}\sup_{\|f\|_1=1}\sup_{\lambda>0}\lambda\cdot \mu\{x:
\sup_n\frac{\bN_n(f)(x)}{n}>\lambda\}\leq C<\infty.$$ This would
contradict Theorem \ref{th1}.
\end{proof}

\subsection{Proof of Theorem \ref{th3}}\label{sec32}
\begin{proof}
  Theorem \ref{th3} also follows from Theorem \ref{th1}.
We can argue also by
  contradiction.
  If we had the validity of the Return Times for Pairs
property
for $(L^1,
  L^1)$ spaces then we would have the convergence in the universal
  sense of the averages
  $$\frac{1}{N}\sum_{n=1}^N f(T^nx)\cdot g(S^ny)= \frac{\sigma _N}{N}$$
  for $g\in L^1(\nu)$.
   This would imply the convergence to zero of
   $$\frac{\sigma _N}{N}-
\frac{\sigma _{N-1}}{N-1}= \frac{f(T^Nx)g(S^Ny)}{N} +
   \frac{\sigma
_{N-1}}{N-1}\cdot \frac{N-1}{N}- \frac{\sigma _{N-1}}{N-1}.$$
This in turn would give the validity of the Return Time for the
Tail property in $L^1$, but this was disproved in Theorem \ref{th2}.
\end{proof}

\section{The counting problem and Birkhoff's theorem}\label{sec4}
  Theorem \ref{th1} also helps to
 refine connections between Birkhoff's pointwise
  ergodic theorem and the counting problem. It provides an
  example of a maximal operator which is of restricted weak type (1,1)
  but does not satisfy a weak type (1, 1) inequality. However, this operator coincides with
  the one sided Hardy--Littlewood maximal function on characteristic functions of measurable sets. Let us see
  how and why.
  \vskip1ex
  One way to prove Birkhoff's pointwise ergodic theorem is via the
  maximal inequality
  $$\mu\left \{x: \sup_N\frac{1}{N}\sum_{n=1}^N |f|(T^nx)>\lambda
\right\}\leq
  \frac{1}{\lambda}\|f\|_1.$$
  It turns out (see \cite{[BS]} for instance)
 that this maximal inequality is equivalent to the
  weak type (1,1) inequality for the Hardy--Littlewood maximal function
 on $\T$, the unit circle, that we identify with the interval
 $[-\frac{1}{2}, \frac{1}{2})$,
 $$H(f)(x) = \sup_{t>0}\frac{1}{t}\int_0^t |f(x-y)|dy.$$

 The following maximal function was introduced by the first
 author
$$A(f)(x) = \sup_{\lambda>0}\lambda\cdot  m\left \{0<y<x:
\frac{|f(x-y)|}{y}>\lambda\right\}.$$

 The interest in the operator $A$ lies in the following results

\begin{enumerate}
\item It was used in \cite{[N]} to give the details of the fact
that the return time for the tail in all $L^p$ spaces $1<p\leq
\infty$ is equivalent to the validity of Birkhoff's theorem in all
$L^r$ spaces for $1<r\leq \infty$.  In other words, the finiteness
of $\bN^*(f)(x) = \sup_n\frac{\bN_n(f)(x)}{n}$ shown in
\cite{[A1]} is equivalent to Birkhoff's theorem in $L^p$ for
$1<p\leq \infty$.
 \item If one considers the characteristic function of a
measurable set $B,$ then simple computations show that
\begin{equation}
A(\I_{B})(x)= H(\I_B)(x). \end{equation}
  Thus the operator $A$ satisfies a restricted weak type (1, 1)
 inequality in the sense that we have for all
 $\lambda>0$
 $$m\{x: A(\I_B)(x) >\lambda\}\leq \frac{1}{\lambda}m(B)$$ i.e. a weak
 type (1, 1) inequality for characteristic functions of
 measurable sets. (See also \cite{[SW]} or \cite{[BS]} for instance for more on
 restricted weak type inequalities.)
 \item The operator $A$ can be viewed as a continuous analog of
 the counting function studied in the previous sections. Furthermore,
 we have
 the following lemma.

\begin{lemma}\label{lem8}
Given $p$, $1\leq p \leq \infty$ the following statements are
equivalent
\begin{enumerate}
\item There exists a finite constant $C$ such that for all
$\lambda>0$ and $(a_n)\in l^p(\Z)$
\begin{equation}
\#\left \{i\in \Z: \sup_n\left (\frac{\#\{k>0:
\frac{a_{k+i}}{k}>\frac{1}{n}\}}{n}\right)>\lambda\right\}\leq
\frac{C}{\lambda^p}\|(a_n)\|_p^p.
\end{equation}
 \item There exists a finite constant $C$ such that for all
$f\in L^p(\T)$ and $\lambda>0$ we have
$$m\{x: A(f)(x)>\lambda\} \leq \frac{C}{\lambda^p}\int |f|^p dm.$$
\item We can find a finite constant $C$ such that for all $f\in
L_+^p(\mu)$ for all measure preserving systems $(X, \mathcal{B},
\mu, T)$
$$\mu\left \{x: \sup_n\frac{\bN_n(f)(x)}{n}>\lambda\right\}\leq
\frac{C}{\lambda^p}\int |f|^p d\mu$$
\end{enumerate}
\end{lemma}

\begin{proof}
 The proof uses known methods in ergodic theory such as
 transference or Rohlin's tower lemma.
Details of such computations can be seen in \cite{[N]}.
So we only sketch
 some of them.
  We remark that (a) is equivalent to the following inequality.
  \vskip1ex
  There exists a finite constant $C$ such that
 for all $\lambda>0$, $(a_n)\in l^p(\Z)$, positive integers $K$ and $I$,
 \begin{equation}
\#\left \{i\in [-I,I]: \sup_{n\leq K}\left (\frac{\#\{k>0:
\frac{a_{k+i}}{k} >\frac{1}n\} }{n}\right)>\lambda\right\}\leq
\frac{C}{\lambda^p}\|(a_n)\|_p^p.
\end{equation}

 In order to prove that (a) and (b) are equivalent we use step functions of the form $f = \sum_{j=-I}^{I-1}
 a_j{\I}_{I_j}$ where $a_j\in \R$ and $a_j = 0$ for $|j|>I$. The
 interval  $I_i$ equals the dyadic interval $[\frac{i}{2^I},
 \frac{i+1}{2^I})$.
 \vskip1ex
 To show that (a) and (c) are equivalent we use Rohlin's tower lemma
 where the tower is symmetric and of height $2J +1.$ Rohlin's lemma
 tells us that for any $\epsilon>0$ and $J\in \N$ we can find
 disjoint sets $T^{-i}B$ for $-J\leq i\leq J$, such that the tower
 $\cup_{i= -J}^{J} T^{-i}(B)$ has total measure greater than
 $1-\epsilon$ . We take a function $f= \sum_{i=-J}^{J} a_i{\I}_{T^iB}$
 and note that
\begin{align*}
\frac{\bN_n(f)(x)}{n}&= \frac{ \# \{k:\frac{f(T^kx)}{k}>
 \frac{1}{n}\}}{n}\\
 &\geq \sum_{i=-J}^J \I_{T^{i}B}(x)\frac{\#\{k\leq J-|i|:
 \frac{a_{k+i}}{k}> \frac{1}{n}\}}{n}.
\end{align*}
 Thus,
 the inequality
$$\mu\left \{x: \sup_n\frac{\bN_n(f)(x)}{n}>\lambda\right\}\leq
\frac{C}{\lambda^p}\int |f|^p d\mu$$ implies
\begin{align}\label{eq26}
\sum_{i=-J}^{J}\mu&\left \{x\in T^iB: \sup_n\left (\frac{\#\{k\leq J-|i|:
 \frac{a_{k+i}}{k}> \frac{1}{n}\}}{n}\right)>\lambda\right\}\\
\nonumber
&
\leq \frac{C}{\lambda^p}\mu(B)\sum_{i=-J}^J |a_i|^p.
\end{align}
 As \eqref{eq26} equals
 $$\mu(B)\cdot \#\left \{-J\leq i\leq J:\sup_n\left (\frac{\#\{k\leq J-|i|:
 \frac{a_{k+i}}{k}> \frac{1}{n}\}}{n}\right)>\lambda\right\}$$
 we have
 \[
\begin{aligned}
 &\#\left \{i\in \Z: \sup_{0<n\leq
 K}\left (\frac{\#\{k>0:\frac{a_{k+i}}{k}\geq \frac{1}{n}\}}{n}
\right)>\lambda \right\} \\
 &\le \lim_J\#\left \{-J\leq i\leq J: \sup_{0<n\leq K}
\left (\frac{\#\{k\leq J-|i|:
 \frac{a_{k+i}}{k}\geq \frac{1}{n}\}}{n}\right )>\lambda\right\}.
 \end{aligned}
 \]
\end{proof}
\end{enumerate}
 So Theorem \ref{th1} gives us the following contribution to the
 problem of characterizing operators for which a
 restricted weak type (1,1) inequality implies a weak type (1,1)
 inequality. (See \cite{[BS]} for more on this problem.)
 The operator $A$ does not satisfy a weak type
 (1,1) inequality. It is shown in \cite{[BS]} that if an operator is generated by
 convolutions, then a restricted weak type (1,1)
 inequality implies a weak type (1,1) inequality. Such is
 the case of the Hilbert transform and the Hardy--Littlewood
 maximal function.\\
  Next we list some of the properties
of the operator $A$.
 \begin{theorem}\label{th9}
The operator A defined on $\T$ by the formula
$$A(f)(x) = \sup_{\lambda>0} \lambda \cdot m \left \{0<y<x:
\frac{|f(x-y)|}{y}>\lambda\right\}$$
 has the following properties
\begin{enumerate}
\item It coincides with the one sided Hardy--Littlewood maximal
function when $f$ is the characteristic function of a measurable
set on $\T$ hence it satisfies a restricted weak type (1,1)
inequality. \item It maps functions in $L^p$ to functions in weak
$L^p$. \item There exists a positive function $f\in L^1(\T)$ such
that $A(f)(x)\nless \infty $ for a.e. x in $\T$.
\end{enumerate}
\end{theorem}
\begin{proof}
Statements (1) and (2) follow from Lemma \ref{lem8}.

The last statement is a consequence of Theorem \ref{th1}. The
arguments developped in \cite {[A1]} (cf. Theorem 4) indicate that
if we had $A(f)(x) <\infty $ for a.e. $x$ then we would have a
weak type (1,1) inequality for $A$. By Lemma \ref{lem8} this would
imply a weak type (1,1) inequality for
$\sup_n\frac{\bN_n(f)(x)}{n}$, a conclusion that we disproved in
Theorem \ref{th1}.
\end{proof}

\end{document}